\documentclass[a4paper, 12 pt]{article}

\usepackage[T2A]{fontenc}
\usepackage[cp1251]{inputenc}
\usepackage[english]{babel}
\usepackage[tbtags]{amsmath}
\usepackage{amsfonts,amssymb}

\usepackage{mathrsfs}

\usepackage{geometry} 
\begin{document}
\title{Composition polynomials of  RNA matrix and $B$-composition  polynomials of Riordan pseudo-involution}
 \author{E. Burlachenko}
 \date{}

 \maketitle
\begin{abstract}
Let $\left( g\left( x \right),xg\left( x \right) \right)$ be a Riordan matrix from the Bell subgroup. We denote ${{\left( g\left( x \right),xg\left( x \right) \right)}^{\varphi }}=\left( {{g}^{\left( \varphi  \right)}}\left( x \right),x{{g}^{\left( \varphi  \right)}}\left( x \right) \right)$, where a  matrix power is defined in the standard way. The polynomials ${{c}_{n}}\left( x \right)$  such that ${{g}^{\left( \varphi  \right)}}\left( x \right)=\sum\nolimits_{n=0}^{\infty }{{{c}_{n}}}\left( \varphi  \right){{x}^{n}}$  will be called composition polynomials. We consider the composition polynomials  of the RNA matrix. The construction associated with these polynomials allows the following generalization. If the matrix $\left( g\left( x \right),xg\left( x \right) \right)$ is a pseudo-involution, then there exists  a numerical sequence ($B$-sequence) with the generating function $B\left( x \right)$ such that $g\left( x \right)=1+xg\left( x \right)B\left( {{x}^{2}}g\left( x \right) \right)$. We derive the formula expressing coefficients of the series ${{g}^{\varphi }}\left( x \right)$ in terms of  coefficients of the series $B\left( x \right)$.  The matrix whose $B$-sequence has the generating function $\varphi B\left( x \right)$ will be denoted by $\left( {{g}^{\left[ \varphi  \right]}}\left( x \right),x{{g}^{\left[ \varphi  \right]}}\left( x \right) \right)$. The polynomials  ${{u}_{n}}\left( x \right)$ such that ${{g}^{\left[ \varphi  \right]}}\left( x \right)=\sum\nolimits_{n=0}^{\infty }{{{u}_{n}}}\left( \varphi  \right){{x}^{n}}$  will be called $B$-composition polynomials. Coefficients of these polynomials are expressed in terms of the $B$-sequence. We show that matrices whose rows correspond to the $B$-composition polynomials are connected with exponential Riordan  matrices of  the Lagrange subgroup in a certain way. The cases  $B\left( x \right)={{\left( 1-x \right)}^{-1}}$ (RNA matrix), $B\left( x \right)=1+x$, $B\left( x \right)=C\left( x \right)$, where $C\left( x \right)$ is the Catalan series, are considered in detail. 
\end{abstract}

\section{Introduction}

Matrices that we will consider correspond to operators in the space of formal power series over the field of real or complex  numbers. Based on this, we associate the rows  and  columns of matrices with the generating functions of their elements, i.e., formal power series. Thus, the expression $Aa\left( x \right)=b\left( x \right)$ means that the column vector multiplied by the matrix $A$ has the generating function $a\left( x \right)=\sum\nolimits_{n=0}^{\infty }{{{a}_{n}}}{{x}^{n}}$, resultant column vector has the generating function $b\left( x \right)=\sum\nolimits_{n=0}^{\infty }{{{b}_{n}}{{x}^{n}}}$. The $n$th coefficient of the series $a\left( x \right)$ will be denoted by $\left[ {{x}^{n}} \right]a\left( x \right)$; the $\left( n,m \right)$th element of the matrix $A$, generating functions of the $n$th row, $n$th descending diagonal, $n$th ascending diagonal and the $n$th column of the matrix $A$ will be denoted  respectively by
$${{\left( A \right)}_{n,m}},\qquad \left[ n,\to  \right]A, \qquad\left[ n,\searrow  \right]A,   \qquad[n,\nearrow ]A,   \qquad A{{x}^{n}}.$$

The infinite lower triangular matrix $\left( f\left( x \right),g\left( x \right) \right)$, the $n$th column of which has the generating function $f\left( x \right){{g}^{n}}\left( x \right)$, ${{g}_{0}}=0$, is called Riordan matrix (Riordan array) [2,12]. It is the product of two matrices that correspond to the operations of multiplication and composition of series:
$$\left( f\left( x \right),g\left( x \right) \right)=\left( f\left( x \right),x \right)\left( 1,g\left( x \right) \right),$$
$$\left( f\left( x \right),x \right)a\left( x \right)=f\left( x \right)a\left( x \right), \quad\left( 1,g\left( x \right) \right)a\left( x \right)=a\left( g\left( x \right) \right),$$
$$\left( f\left( x \right),g\left( x \right) \right)\left( b\left( x \right),a\left( x \right) \right)=\left( f\left( x \right)b\left( g\left( x \right) \right),a\left( g\left( x \right) \right) \right).$$

Matrices $\left( f\left( x \right),g\left( x \right) \right)$, ${{f}_{0}}\ne 0$, ${{g}_{1}}\ne 0$, or in the more convenient notation for our topic $\left( f\left( x \right),xg\left( x \right) \right)$, ${{f}_{0}}\ne 0$, ${{g}_{0}}\ne 0$, form a group called the Riordan group. Matrices of the form $\left( f\left( x \right),x \right)$ form a subgroup called the Appell subgroup, matrices of the form form $\left( 1,xg\left( x \right) \right)$ form a subgroup called the Lagrange subgroup, or associated subgroup. A subgroup of the matrices $\left( g\left( x \right),xg\left( x \right) \right)$ isomorphic to the Lagrange subgroup is called the Bell subgroup. 

Еlements of Riordan matrix will be denoted by ${{d}_{n,m}}$. If  ${{d}_{n,n}}\ne 0$, then there exist a unique numerical sequence $A={{\left\{ {{a}_{n}} \right\}}_{n\ge 0}}$ called $A$-sequence such that 
$${{d}_{n+1,m+1}}=\sum\limits_{i=0}^{\infty }{{{a}_{i}}{{d}_{n,m+i}}},$$
[10,17]. Let $A\left( x \right)$ be the generating function of the $A$-sequence of the matrix $\left( f\left( x \right),xg\left( x \right) \right)$. Then
$$f\left( x \right){{g}^{m+1}}\left( x \right)=f\left( x \right){{g}^{m}}\left( x \right)A\left( xg\left( x \right) \right),$$
$$g\left( x \right)=A\left( xg\left( x \right) \right),  \qquad{{\left( 1,xg\left( x \right) \right)}^{-1}}=\left( 1,x{{A}^{-1}}\left( x \right) \right).$$
For example, $A\left( x \right)=1+{{a}_{1}}x+{{a}_{2}}{{x}^{2}}$, 
$$g\left( x \right)=1+{{a}_{1}}xg\left( x \right)+{{a}_{2}}{{x}^{2}}{{g}^{2}}\left( x \right)=\frac{1-{{a}_{1}}x-\sqrt{{{\left( 1-{{a}_{1}}x \right)}^{2}}-4{{a}_{2}}{{x}^{2}}}}{2{{a}_{2}}{{x}^{2}}}$$.

Matrices of the form 
$$\left( a\left( x \right),x \right)+\left( xb\left( x \right),x \right)D,$$
where $D$ is the matrix of the differentiation operator: $D{{x}^{n}}=n{{x}^{n-1}}$, and $a\left( x \right)$, $b\left( x \right)$ are arbitrary series, form the Lie algebra of  the Riordan group [1,19]. 

The matrices 
$${{\left| {{e}^{x}} \right|}^{-1}}\left( f\left( x \right),xg\left( x \right) \right)\left| {{e}^{x}} \right|={{\left( f\left( x \right),xg\left( x \right) \right)}_{E}},$$
where $\left| {{e}^{x}} \right|$ is the diagonal matrix: $\left| {{e}^{x}} \right|{{x}^{n}}={{{x}^{n}}}/{n!}\;$, are called exponential Riordan matrices [21,22]. Denote $\left[ n,\to  \right]{{\left( f\left( x \right),xg\left( x \right) \right)}_{E}}={{s}_{n}}\left( x \right)$. Then
$$\sum\limits_{n=0}^{\infty }{\frac{{{s}_{n}}\left( \varphi  \right)}{n!}{{x}^{n}}}=f\left( x \right)\exp \left( \varphi xg\left( x \right) \right).$$
In a general case (but for ${{f}_{0}},{{g}_{0}}\ne 0$) the sequence of polynomials ${{s}_{n}}\left( x \right)$ is called the Scheffer sequence; in the case $g\left( x \right)=1$ –  the Appell sequence, in the case $f\left( x \right)=1$ –  the binomial sequence. The matrix $P$ whose power is defined by the identity
$${{P}^{\varphi }}=\left( \frac{1}{1-\varphi x},\frac{x}{1-\varphi x} \right)={{\left( {{e}^{\varphi x}},x \right)}_{E}}$$
is called the Pascal matrix. 

The Riordan matrix $\left( f\left( x \right),xg\left( x \right) \right)$, ${{g}_{0}}=\pm 1$, having the property
$${{\left( f\left( x \right),xg\left( x \right) \right)}^{-1}}=\left( 1,-x \right)\left( f\left( x \right),xg\left( x \right) \right)\left( 1,-x \right)=\left( f\left( -x \right),xg\left( -x \right) \right)$$
is called pseudo-involution in the Riordan group [3,5,6,7,8,10,11,16,18]. The case ${{g}_{0}}=-1$ corresponds to the matrices $\left( 1,-x \right)$, $\left( -1,-x \right)$, which are both involutions and pseudo-involutions. Obviously, if the Riordan matrix $A$ is a pseudo-involution, then the matrices $\left( 1,-x \right)A$, $A\left( 1,-x \right)$ are involutions. A example of  pseudo-involution is the power of Pascal matrix. The power of pseudo-involution, as well as any polyindromic product of pseudo-involutions (i.e.,  product that does not change when its members are rearranged in the reverse order) is also a pseudo-involution [11].

If the matrix $\left( f\left( x \right),xg\left( x \right) \right)$, ${{g}_{0}}=1$, is a pseudo-involution, then there exist a unique numerical sequence $B={{\left\{ {{b}_{n}} \right\}}_{n\ge 0}}$ such that 
$${{d}_{n+1,m}}={{d}_{n,m-1}}+\sum\limits_{i=0}^{\infty }{{{b}_{i}}{{d}_{n-i,m+i}}},  \qquad{{d}_{n,-1}}=0,$$
[8,18]. Let $B\left( x \right)$ be the generating function of this sequence. Then
$$f\left( x \right){{g}^{m}}\left( x \right)=f\left( x \right){{g}^{m-1}}\left( x \right)+xf\left( x \right){{g}^{m}}\left( x \right)B\left( {{x}^{2}}g\left( x \right) \right),$$
$$g\left( x \right)=1+xg\left( x \right)B\left( {{x}^{2}}g\left( x \right) \right).$$
For example, $B\left( x \right)={{b}_{0}}+{{b}_{1}}x$,  
$$g\left( x \right)=1+{{b}_{0}}xg\left( x \right)+{{b}_{1}}{{x}^{3}}{{g}^{2}}\left( x \right)=\frac{1-{{b}_{0}}x-\sqrt{{{\left( 1-{{b}_{0}}x \right)}^{2}}-4{{b}_{1}}{{x}^{3}}}}{2{{b}_{1}}{{x}^{3}}}.$$
The sequence $B$ is called $B$-sequence of the matrix $\left( f\left( x \right),xg\left( x \right) \right)$ (in [8] this sequence is called $\Delta $-sequence). The generating function of this sequence will be called  $B$-function  of the matrix $\left( f\left( x \right),xg\left( x \right) \right)$. 

Consider the following construction for the Bell subgroup matrices $\left( g\left( x \right),xg\left( x \right) \right)$, ${{g}_{0}}=1$. Denote
$${{\left( g\left( x \right),xg\left( x \right) \right)}^{\varphi }}=\sum\limits_{n=0}^{\infty }{\left( \begin{matrix}
   \varphi   \\
   n  \\
\end{matrix} \right)}{{\left( \left( g\left( x \right),xg\left( x \right) \right)-I \right)}^{n}},$$
$$\log \left( g\left( x \right),xg\left( x \right) \right)=\sum\limits_{n=1}^{\infty }{\frac{{{\left( -1 \right)}^{n-1}}}{n}}{{\left( \left( g\left( x \right),xg\left( x \right) \right)-I \right)}^{n}}.$$
where $I=\left( 1,x \right)$. Then
$${{\left( g\left( x \right),xg\left( x \right) \right)}^{\varphi }}=\sum\limits_{n=0}^{\infty }{\frac{{{\varphi }^{n}}}{n!}}{{\left( \log \left( g\left( x \right),xg\left( x \right) \right) \right)}^{n}}.$$
We will build the matrix $L\left( g\left( x \right) \right)$ by the rule 
$$L\left( g\left( x \right) \right){{x}^{n}}=\left( {1}/{n!}\; \right){{\left( \log \left( g\left( x \right),xg\left( x \right) \right) \right)}^{n}}{{x}^{0}}.$$
Denote ${{\left( g\left( x \right),xg\left( x \right) \right)}^{\varphi }}=\left( {{g}^{\left( \varphi  \right)}}\left( x \right),x{{g}^{\left( \varphi  \right)}}\left( x \right) \right)$, $\left[ n,\to  \right]L\left( g\left( x \right) \right)={{c}_{n}}\left( x \right)$. Then ${{g}^{\left( \varphi  \right)}}\left( x \right)=\sum\nolimits_{n=0}^{\infty }{{{c}_{n}}}\left( \varphi  \right){{x}^{n}}$. The polynomials ${{c}_{n}}\left( x \right)$ will be called composition polynomials. (If in this construction we replace the Bell subgroup matrices with the Appell subgroup matrices, then we get ${{g}^{\left( \varphi  \right)}}\left( x \right)={{g}^{\varphi }}\left( x \right)$, $L\left( g\left( x \right) \right)=\left( 1,\log g\left( x \right) \right)\left| {{e}^{x}} \right|$; in this case the polynomials ${{c}_{n}}\left( x \right)$ are called convolution polynomials [13]). 

The composition polynomials together with the convolution polynomials, irrespective of Riordan matrices, were considered in [14]. With regard to our construction, we note that 
$$\log \left( g\left( x \right),xg\left( x \right) \right)=\left( b\left( x \right),x \right){{D}^{T}},$$
where ${{D}^{T}}=\left( x,x \right)D\left( x,x \right)={{\left( x,x \right)}_{E}}$, and the series $b\left( x \right)$ satisfies the conditions ${{b}_{0}}={{g}_{1}}$, ${{g}^{2}}\left( x \right)b\left( xg\left( x \right) \right)=b\left( x \right){{\left( xg\left( x \right) \right)}^{\prime }}$. Thus, 
$$L\left( g\left( x \right) \right){{x}^{0}}=1, \qquad L\left( g\left( x \right) \right){{x}^{n}}=\left( {1}/{n}\; \right)b\left( x \right){{D}^{T}}L\left( g\left( x \right) \right){{x}^{n-1}}.$$
Coefficients of the polynomials ${{c}_{n}}\left( \varphi  \right)=\left[ {{x}^{n}} \right]{{g}^{\left( \varphi  \right)}}\left( x \right)$ are expressed in terms of  coefficients of the series $b\left( x \right)$ by the formula
$${{c}_{n}}\left( \varphi  \right)=\sum\limits_{m=0}^{n}{\frac{{{\varphi }^{m}}}{m!}}\sum\limits_{n,m}{{{b}_{{{i}_{1}}-1}}}{{b}_{{{i}_{2}}-1}}...{{b}_{{{i}_{m}}-1}}\left( 1+{{i}_{1}} \right)\left( 1+{{i}_{1}}+{{i}_{2}} \right)...\left( 1+{{i}_{1}}+{{i}_{2}}+...+{{i}_{m-1}} \right),$$
where the summation of the coefficient of ${{{\varphi }^{m}}}/{m!}\;$ is over all compositions $n={{i}_{1}}+{{i}_{2}}+...+{{i}_{m}}$ . A generalization is the formula 
$${{c}_{n}}\left( \beta ,\varphi  \right)=\left[ {{x}^{n}} \right]{{\left( {{g}^{\left( \varphi  \right)}}\left( x \right) \right)}^{\beta }}=$$
$$=\sum\limits_{m=0}^{n}{\frac{{{\varphi }^{m}}}{m!}}\sum\limits_{n,m}{{{b}_{{{i}_{1}}-1}}}{{b}_{{{i}_{2}}-1}}...{{b}_{{{i}_{m}}-1}}\beta \left( \beta +{{i}_{1}} \right)\left( \beta +{{i}_{1}}+{{i}_{2}} \right)...\left( \beta +{{i}_{1}}+{{i}_{2}}+...+{{i}_{m-1}} \right).$$

Note that if the matrix $\left( g\left( x \right),xg\left( x \right) \right)$ is a pseudo-involution, i.e., ${{g}^{\left( -1 \right)}}\left( x \right)=g\left( -x \right)$, then the polynomial ${{c}_{2n}}\left( x \right)$ is  an even function, the polynomial ${{c}_{2n+1}}\left( x \right)$ is an odd function. Thus, the matrix $\left( b\left( x \right),x \right){{D}^{T}}$ is the logarithm of  pseudo-involution if $b\left( x \right)$ is an even function. \\
{\bfseries Example 1.1.}
$$g\left( x \right)={{\left( 1-x \right)}^{-1}}, \qquad\left( g\left( x \right),xg\left( x \right) \right)={{\left( {{e}^{x}},x \right)}_{E}},$$
$$\log \left( g\left( x \right),xg\left( x \right) \right)={{\left( x,x \right)}_{E}}, \qquad L\left( g\left( x \right) \right)=\left( 1,x \right), \qquad{{c}_{n}}\left( x \right)={{x}^{n}}.$$
{\bfseries Example 1.2.}
$${{\left( g\left( x \right),xg\left( x \right) \right)}^{\varphi }}=A{{P}^{\varphi }}{{A}^{-1}}=$$
$$=\left( \frac{1}{\sqrt{1-{{x}^{2}}}},\frac{x}{\sqrt{1-{{x}^{2}}}} \right)\left( \frac{1}{1-\varphi x},\frac{x}{1-\varphi x} \right)\left( \frac{1}{\sqrt{1+{{x}^{2}}}},\frac{x}{\sqrt{1+{{x}^{2}}}} \right),$$
$${{g}^{\left( \varphi  \right)}}\left( x \right)=\frac{1}{\sqrt{1-2\varphi x\sqrt{1-{{x}^{2}}}+{{\varphi }^{2}}{{x}^{2}}}}.$$
Since the matrix ${{P}^{\varphi }}$ is a pseudo-involution and $A\left( 1,-x \right){{A}^{-1}}=\left( 1,-x \right)$, then the matrix ${{\left( g\left( x \right),xg\left( x \right) \right)}^{\varphi }}$ is also a pseudo-involution. Since ${{g}^{\left( \varphi  \right)}}\left( x \right)=\sum\nolimits_{n=0}^{\infty }{{{\varphi }^{n}}{{x}^{n}}}{{P}_{n}}\left( \sqrt{1-{{x}^{2}}} \right)$, where ${{P}_{n}}\left( x \right)$ are the Legendre polynomials, then in this case $b\left( x \right)=\sqrt{1-{{x}^{2}}}$, $L\left( g\left( x \right) \right){{x}^{n}}={{x}^{n}}{{P}_{n}}\left( \sqrt{1-{{x}^{2}}} \right)$.

In Section 2, we associate the $B$-sequence of the matrix $\left( 1,xg\left( x \right) \right)$ with the $A$-sequence of the matrix $\left( 1,x\sqrt{g\left( x \right)} \right)$ and consider examples of such connection. In Section 3, for the  matrix $\left( 1,xg\left( x \right) \right)$ with the $B$-function $B\left( x \right)$, we express coefficients of the series ${{g}^{\varphi }}\left( x \right)$ in terms of coefficients of the series $B\left( x \right)$. In Section 4, we consider the composition polynomials of the RNA matrix which, after the Pascal matrix, is the most famous example of a pseudo-involution in the Riordan group. The construction associated with these polynomials allows the generalization which we introduce in Section 3. The Bell subgroup matrix whose $B$-sequence has the generating function $\varphi B\left( x \right)$ will be denoted by $\left( {{g}^{\left[ \varphi  \right]}}\left( x \right),x{{g}^{\left[ \varphi  \right]}}\left( x \right) \right)$. The polynomials  ${{u}_{n}}\left( x \right)$ such that ${{g}^{\left[ \varphi  \right]}}\left( x \right)=\sum\nolimits_{n=0}^{\infty }{{{u}_{n}}}\left( \varphi  \right){{x}^{n}}$  will be called $B$-composition polynomials. Coefficients of these polynomials are expressed in terms of the $B$-sequence. The matrix whose rows correspond to the $B$-composition polynomials will be called $B$-composition matrix. In Section 6, Section 7, we will build the $B$-composition matrices for the cases $B\left( x \right)=1+x$, $B\left( x \right)=C\left( x \right)$, where $C\left( x \right)$ is the Catalan series. Both cases, as well as the case $B\left( x \right)={{\left( 1-x \right)}^{-1}}$(RNA matrix), are related to the Narayana polynomials in a certain way. In Section 8, we prove the simple but unexpected theorem on the connection of the $B$-composition matrices with exponential Riordan matrices of the Lagrange subgroup. Using this connection, in Section 9 we introduce a $B$-composition-convolution polynomials ${{u}_{n}}\left( \beta ,x \right)$ such that   $\left[ {{x}^{n}} \right]{{\left( {{g}^{\left[ \varphi  \right]}}\left( x \right) \right)}^{\beta }}={{u}_{n}}\left( \beta ,\varphi  \right)$
\section{Some examples}
{\bfseries Lemma 2.1.} \emph{If the matrix $\left( 1,xg\left( x \right) \right)$, $g\left( x \right)\ne -1$,  is a pseudo-involution, i.e., 
$${{\left( 1,xg\left( x \right) \right)}^{-1}}=\left( 1,-x \right)\left( 1,xg\left( x \right) \right)\left( 1,-x \right)=\left( 1,xg\left( -x \right) \right),$$
then it can be represented in the form
$$\left( 1,xg\left( x \right) \right)=\left( 1,x\sqrt{g\left( x \right)} \right)\left( 1,xh\left( x \right) \right),$$
where $h\left( x \right)$ is the generating function of the $A$-sequence of the matrix $\left( 1,x\sqrt{g\left( x \right)} \right)$ such that}
$$h\left( -x \right)={{h}^{-1}}\left( x \right),   \qquad h\left( x \right)=s\left( x \right)+\sqrt{{{s}^{2}}\left( x \right)+1}, \quad{{s}_{2n}}=0.$$
{\bfseries Proof.} If
$${{\left( 1,x\sqrt{g\left( x \right)} \right)}^{-1}}=\left( 1,x{{h}^{-1}}\left( x \right) \right), \qquad{{\left( 1,xh\left( x \right) \right)}^{-1}}=\left( 1,x\sqrt{c\left( x \right)} \right),$$
then
$$\left( 1,x\sqrt{g\left( x \right)} \right)\left( 1,xh\left( x \right) \right)=\left( 1,xg\left( x \right) \right),$$ 
$$\left( 1,x\sqrt{c\left( x \right)} \right)\left( 1,x{{h}^{-1}}\left( x \right) \right)=\left( 1,xc\left( x \right) \right),  \qquad{{\left( 1,xg\left( x \right) \right)}^{-1}}=\left( 1,xc\left( x \right) \right).$$
It follows from the condition $c\left( x \right)=g\left( -x \right)$ that ${{h}^{-1}}\left( x \right)=h\left( -x \right)$.   \qquad    $\square $ \\
{\bfseries Example 2.1.}
$$\left( 1,\frac{x}{1-2\varphi x} \right)=\left( 1,\frac{x}{\sqrt{1-2\varphi x}} \right)\left( 1,x\left( \varphi x+\sqrt{{{\varphi }^{2}}{{x}^{2}}+1} \right) \right).$$
{\bfseries Example 2.2.}
$$\left( 1,x\sum\limits_{n=0}^{\infty }{\frac{2{{\left( 2+n \right)}^{n-1}}}{n!}{{\varphi }^{n}}{{x}^{n}}} \right)=\left( 1,x\sum\limits_{n=0}^{\infty }{\frac{{{\left( 1+n \right)}^{n-1}}}{n!}{{\varphi }^{n}}{{x}^{n}}} \right)\left( 1,x{{e}^{\varphi x}} \right),$$
where
$$x\sum\limits_{n=0}^{\infty }{\frac{{{\left( 1+n \right)}^{n-1}}}{n!}{{\varphi }^{n}}{{x}^{n}}=\ln \left( \sum\limits_{n=0}^{\infty }{\frac{{{\left( 1+\varphi n \right)}^{n-1}}}{n!}{{x}^{n}}} \right)}=x{{\left( \sum\limits_{n=0}^{\infty }{\frac{{{\left( 1+\varphi n \right)}^{n-1}}}{n!}{{x}^{n}}} \right)}^{\varphi }},$$
$$\sum\limits_{n=0}^{\infty }{\frac{2{{\left( 2+n \right)}^{n-1}}}{n!}{{\varphi }^{n}}{{x}^{n}}={{\left( \sum\limits_{n=0}^{\infty }{\frac{{{\left( 1+\varphi n \right)}^{n-1}}}{n!}{{x}^{n}}} \right)}^{2\varphi }}}.$$
{\bfseries Example 2.3.}
$$\left( 1,\frac{1-4\varphi x+{{\varphi }^{2}}{{x}^{2}}-\sqrt{{{\left( 1-4\varphi x+{{\varphi }^{2}}{{x}^{2}} \right)}^{2}}-4{{\varphi }^{2}}{{x}^{2}}}}{2{{\varphi }^{2}}x} \right)=$$
$$=\left( 1,\frac{1-\varphi x-\sqrt{{{\left( 1-\varphi x \right)}^{2}}-4\varphi x}}{2\varphi } \right)\left( 1,x\frac{1+\varphi x}{1-\varphi x} \right),$$
$$\frac{1+\varphi x}{1-\varphi x}=\frac{2\varphi x}{1-{{\varphi }^{2}}{{x}^{2}}}+\sqrt{{{\left( \frac{2\varphi x}{1-{{\varphi }^{2}}{{x}^{2}}} \right)}^{2}}+1}.$$
{\bfseries Theorem 2.2.} \emph{If $B\left( x \right)$ is the $B$-function of the matrix $\left( 1,xg\left( x \right) \right)$, then in the notation of  Lemma 2.1
$$xB\left( {{x}^{2}} \right)=2s\left( x \right).$$}
{\bfseries Proof.} Since ${{h}^{2}}\left( x \right)=1+2s\left( x \right)h\left( x \right)$, then
$$g\left( x \right)=\left( 1,x\sqrt{g\left( x \right)} \right)\left( 1+2s\left( x \right)h\left( x \right) \right)=1+xg\left( x \right)\tilde{s}\left( x\sqrt{g\left( x \right)} \right)=$$
$$=1+xg\left( x \right)B\left( {{x}^{2}}g\left( x \right) \right),  \qquad\tilde{s}\left( x \right)=\frac{2s\left( x \right)}{x}.\qquad \square $$
{\bfseries Example 2.4.}
 The paper [18] contains the interesting fact that if
$$g\left( x \right)=\sum\limits_{n=0}^{\infty }{\frac{2m+1}{2m+1+\left( m+1 \right)n}}\left( \begin{matrix}
   2m+1+\left( m+1 \right)n  \\
   n  \\
\end{matrix} \right){{x}^{n}},$$
then $B$-sequence of the matrix $\left( 1,xg\left( x \right) \right)$ coincides with the $m$th row of the matrix 
$$\left( \frac{1+x}{{{\left( 1-x \right)}^{2}}},\frac{x}{{{\left( 1-x \right)}^{2}}} \right)=\left( \begin{matrix}
   1 & 0 & 0 & 0 & \cdots   \\
   3 & 1 & 0 & 0 & \cdots   \\
   5 & 5 & 1 & 0 & \cdots   \\
   7 & 14 & 7 & 1 & \cdots   \\
   \vdots  & \vdots  & \vdots  & \vdots  & \ddots   \\
\end{matrix} \right).$$
This is a consequence of the fact that in this case
$$h\left( x \right)={{\left( \frac{x+\sqrt{{{x}^{2}}+4}}{2} \right)}^{2m+1}},\quad
{{\left( \frac{x+\sqrt{{{x}^{2}}+4}}{2} \right)}^{n}}=\frac{{{c}_{n}}\left( x \right)+{{s}_{n-1}}\left( x \right)\sqrt{{{x}^{2}}+4}}{2},$$
$${{s}_{2m}}\left( x \right)\sqrt{{{x}^{2}}+4}=\sqrt{c_{2m+1}^{2}\left( x \right)+4},  \qquad{{c}_{2m}}\left( x \right)=\sqrt{s_{2m-1}^{2}\left( x \right)\left( {{x}^{2}}+4 \right)+4},$$
where the polynomial ${{c}_{n}}\left( x \right)$, $n>0$, corresponds to the $n$th row of the matrix
$$\left( \frac{1+{{x}^{2}}}{1-{{x}^{2}}},\frac{x}{1-{{x}^{2}}} \right)=\left( \begin{matrix}
   1 & 0 & 0 & 0 & 0 & 0 & \cdots   \\
   0 & 1 & 0 & 0 & 0 & 0 & \cdots   \\
   2 & 0 & 1 & 0 & 0 & 0 & \cdots   \\
   0 & 3 & 0 & 1 & 0 & 0 & \cdots   \\
   2 & 0 & 4 & 0 & 1 & 0 & \cdots   \\
   0 & 5 & 0 & 5 & 0 & 1 & \cdots   \\
   \vdots  & \vdots  & \vdots  & \vdots  & \vdots  & \vdots  & \ddots   \\
\end{matrix} \right),$$
the polynomial ${{s}_{n}}\left( x \right)$ corresponds to the $n$th row of the matrix
$$\left( \frac{1}{1-{{x}^{2}}},\frac{x}{1-{{x}^{2}}} \right)=\left( \begin{matrix}
   1 & 0 & 0 & 0 & 0 & 0 & \cdots   \\
   0 & 1 & 0 & 0 & 0 & 0 & \cdots   \\
   1 & 0 & 1 & 0 & 0 & 0 & \cdots   \\
   0 & 2 & 0 & 1 & 0 & 0 & \cdots   \\
   1 & 0 & 3 & 0 & 1 & 0 & \cdots   \\
   0 & 3 & 0 & 4 & 0 & 1 & \cdots   \\
   \vdots  & \vdots  & \vdots  & \vdots  & \vdots  & \vdots  & \ddots   \\
\end{matrix} \right).$$
\section{$B$-expansion}
The generating function of the $n$th row of the matrix $\left( 1,f\left( x \right) \right)$, ${{f}_{0}}=0$, $n>0$, has the form
$$\sum\limits_{n}{\frac{q!{{x}^{q}}}{{{m}_{1}}!{{m}_{2}}!...{{m}_{n}}!}f_{1}^{{{m}_{1}}}f_{2}^{{{m}_{2}}}...f_{n}^{{{m}_{n}}}},   \qquad q=\sum\limits_{i=1}^{n}{{{m}_{i}}}.$$
where the summation is over all partitions $n=\sum\nolimits_{i=1}^{n}{{{m}_{i}}i}$. Then if 
$$g\left( x \right)=a\left( f\left( x \right) \right),  \qquad{{a}^{\varphi }}\left( x \right)=\sum\limits_{n=0}^{\infty }{\frac{{{s}_{n}}\left( \varphi  \right)}{n!}{{x}^{n}}}, \qquad{{s}_{n}}\left( x \right)=\left[ n,\to  \right]{{\left( 1,\log a\left( x \right) \right)}_{E}},$$
then
$$\left[ {{x}^{n}} \right]{{g}^{\varphi }}\left( x \right)=g_{n}^{\left( \varphi  \right)}=\sum\limits_{n}{\frac{{{s}_{q}}\left( \varphi  \right)}{{{m}_{1}}!{{m}_{2}}!...{{m}_{n}}!}f_{1}^{{{m}_{1}}}f_{2}^{{{m}_{2}}}...f_{n}^{{{m}_{n}}}}.$$
A representation of the coefficients $g_{n}^{\left( \varphi  \right)}$ in this form (with the unspoken condition $g_{0}^{\left( \varphi  \right)}=1$) will be called expansion of  binomial type, or binomial expansion. For example, since
$$g\left( x \right)=\left( 1,g\left( x \right)-1 \right)\left( 1+x \right)=\left( 1,\log g\left( x \right) \right){{e}^{x}},$$
then
$$g_{n}^{\left( \varphi  \right)}=\sum\limits_{n}{\frac{{{\left( \varphi  \right)}_{q}}}{{{m}_{1}}!{{m}_{2}}...{{m}_{n}}!}}g_{1}^{{{m}_{1}}}g_{2}^{{{m}_{2}}}...g_{n}^{{{m}_{n}}}=\sum\limits_{n}{\frac{{{\varphi }^{q}}}{{{m}_{1}}!{{m}_{2}}!\text{ }...\text{ }{{m}_{n}}!}}\text{ }l_{1}^{{{m}_{1}}}l_{2}^{{{m}_{2}}}...\text{ }l_{n}^{{{m}_{n}}},$$
where  ${{\left( \varphi  \right)}_{q}}=\varphi \left( \varphi -1 \right)...\left( \varphi -q+1 \right)$, ${{l}_{i}}=\left[ {{x}^{i}} \right]\log g\left( x \right)$. Denote
$$\left[ n,\to  \right]\left( 1,\log g\left( x \right) \right)\left| {{e}^{x}} \right|={{l}_{n}}\left( x \right),  \qquad\left[ n,\to  \right]\left( 1,\log A\left( x \right) \right)\left| {{e}^{x}} \right|={{\tilde{l}}_{n}}\left( x \right).$$
Since ${{\left( 1,xg\left( x \right) \right)}^{-1}}=\left( 1,x{{A}^{-1}}\left( x \right) \right)$, then by the Lagrange inversion theorem
$${{l}_{n}}\left( x \right)=\frac{x}{x+n}{{\tilde{l}}_{n}}\left( x+n \right).$$
Thus,
$${{\tilde{l}}_{n}}\left( x \right)=\sum\limits_{n}^{{}}{\frac{{{\left( x \right)}_{q}}}{{{m}_{1}}!{{m}_{2}}...{{m}_{n}}!}}a_{1}^{{{m}_{1}}}a_{2}^{{{m}_{2}}}...a_{n}^{{{m}_{n}}},  \quad{{l}_{n}}\left( x \right)=\sum\limits_{n}^{{}}{\frac{x{{\left( x+n-1 \right)}_{q-1}}}{{{m}_{1}}!{{m}_{2}}...{{m}_{n}}!}}a_{1}^{{{m}_{1}}}a_{2}^{{{m}_{2}}}...a_{n}^{{{m}_{n}}}.$$
The expansion
$$g_{n}^{\left( \varphi  \right)}=\sum\limits_{n}{\frac{\varphi {{\left( \varphi +n-1 \right)}_{q-1}}}{{{m}_{1}}!{{m}_{2}}...{{m}_{n}}!}}a_{1}^{{{m}_{1}}}a_{2}^{{{m}_{2}}}...a_{n}^{{{m}_{n}}}\eqno (1)$$
is not expansion of binomial type; therefore, we will give a special name to such expansions. The expansion such that 
$$g_{n}^{\left( \varphi  \right)}=\sum\limits_{n}{\frac{\left( {\varphi }/{\beta }\; \right){{s}_{q}}\left( \left( {\varphi }/{\beta }\; \right)+n \right)}{\left( \left( {\varphi }/{\beta }\; \right)+n \right){{m}_{1}}!{{m}_{2}}!...{{m}_{n}}!}f_{1}^{{{m}_{1}}}f_{2}^{{{m}_{2}}}...f_{n}^{{{m}_{n}}}},$$
if
$$\left[ {{x}^{n}} \right]{}_{\left( \beta  \right)}{{A}^{\varphi }}\left( x \right)=\sum\limits_{n}{\frac{{{s}_{q}}\left( \varphi  \right)}{{{m}_{1}}!{{m}_{2}}!...{{m}_{n}}!}f_{1}^{{{m}_{1}}}f_{2}^{{{m}_{2}}}...f_{n}^{{{m}_{n}}}},$$
where $_{\left( \beta  \right)}A\left( x \right)$ is the generating function of the $A$-sequence of the matrix $\left( 1,x{{g}^{\beta }}\left( x \right) \right)$, will be called $A$-binomial expansion.\\
{\bfseries Theorem 3.1.} \emph{
If $B={{\left\{ {{b}_{n}} \right\}}_{n\ge 0}}$ is the $B$-sequence of the matrix $\left( 1,xg\left( x \right) \right)$, then the following formula holds
$$g_{n}^{\left( \varphi  \right)}=\sum\limits_{n}^{{}}{\frac{\varphi {{\left( \varphi +k-1 \right)}_{q-1}}}{{{m}_{0}}!{{m}_{1}}!...{{m}_{p}}!}}b_{0}^{{{m}_{0}}}b_{1}^{{{m}_{1}}}...b_{p}^{{{m}_{p}}}, \eqno (2)$$
$$p=\left\lfloor \frac{n-1}{2} \right\rfloor ,  \qquad k=\sum\limits_{i=0}^{p}{{{m}_{i}}}\left( i+1 \right), \qquad q=\sum\limits_{i=0}^{p}{{{m}_{i}}},$$
where the summation is over all monomials $b_{0}^{{{m}_{0}}}b_{1}^{{{m}_{1}}}...b_{p}^{{{m}_{p}}}$ for which $n=\sum\nolimits_{i=0}^{p}{{{m}_{i}}\left( 2i+1 \right)}$.} 
{\bfseries Proof.}  According to Theorem  2.2.
$${{\left( 1,x\sqrt{g\left( x \right)} \right)}^{-1}}=\left( 1,x{{h}^{-1}}\left( x \right) \right),$$
$$h\left( x \right)=\left( 1,s\left( x \right) \right)\left( x+\sqrt{{{x}^{2}}+1} \right),  \qquad{{s}_{2n}}=0, \qquad{{s}_{2n+1}}={{{b}_{n}}}/{2}\;.$$ 
The binomial expansion of coefficients of  the series ${{h}^{\varphi }}\left( x \right)$ has the form
$$\left[ {{x}^{n}} \right]{{h}^{\varphi }}\left( x \right)=\sum\limits_{n}^{{}}{\frac{{{s}_{q}}\left( \varphi  \right)}{{{m}_{0}}!{{m}_{1}}!...{{m}_{p}}!}}\frac{1}{{{2}^{q}}}b_{0}^{{{m}_{0}}}b_{1}^{{{m}_{1}}}...b_{p}^{{{m}_{p}}},$$ 
where
$$p=\left\lfloor \frac{n-1}{2} \right\rfloor ,  \qquad n=\sum\limits_{i=0}^{p}{{{m}_{i}}}\left( 2i+1 \right), \qquad q=\sum\limits_{i=0}^{p}{{{m}_{i}}},$$
$${{s}_{q}}\left( x \right)=\left[ q,\to  \right]{{\left( 1,\log \left( x+\sqrt{{{x}^{2}}+1} \right) \right)}_{E}}=x\prod\limits_{i=1}^{q-1}{\left( x+q-2i \right)}, \qquad{{s}_{1}}\left( x \right)=x.$$
Since
$$\frac{2\varphi }{2\varphi +n}{{s}_{q}}\left( 2\varphi +n \right)={{2}^{q}}\varphi \prod\limits_{i=1}^{q-1}{\left( \varphi +\frac{q+n}{2}-i \right)}={{2}^{q}}\varphi {{\left( \varphi +k-1 \right)}_{q-1}},$$
where $k=\sum\nolimits_{i=0}^{p}{{{m}_{i}}\left( i+1 \right)}$, then the corresponding $A$-binomial expansion of  coefficients of  the series ${{g}^{\varphi }}\left( x \right)$ has the form (2).        \qquad   $\square $ 

Formula (2) will be called $B$-expansion.\\
{\bfseries Example 3.1.}
Let $B\left( x \right)={{b}_{0}}+{{b}_{r}}{{x}^{r}}$.  The monomial
$$b_{0}^{n-v\left( 2r+1 \right)}b_{r}^{v},\qquad
v=0, 1, … \left\lfloor \frac{n}{2r+1} \right\rfloor ,$$
corresponds to the partition of number $n$ into  parts equal to $1$ and $2r+1$. The coefficient of this monomial in the $B$-expansion is
$$\frac{\varphi {{\left( \varphi +n-v\left( 2r+1 \right)+v\left( r+1 \right)-1 \right)}_{n-v\left( 2r+1 \right)+v-1}}}{\left( n-v\left( 2r+1 \right) \right)!v!}=\frac{\varphi \left( \varphi +n-vr-1 \right)!}{\left( \varphi +vr \right)!\left( n-v\left( 2r+1 \right) \right)!v!}=$$
$$=\left( \begin{matrix}
   \varphi +n-vr-1  \\
   n-v\left( 2r+1 \right)  \\
\end{matrix} \right)\frac{\varphi }{\varphi +rv}\left( \begin{matrix}
   \varphi +\left( r+1 \right)v-1  \\
   v  \\
\end{matrix} \right)=$$
$$=\left( \begin{matrix}
   \varphi +v\left( r+1 \right)+n-v\left( 2r+1 \right)-1  \\
   n-v\left( 2r+1 \right)  \\
\end{matrix} \right)\left[ {{x}^{v}} \right]{{\mathcal{B}}_{r+1}}{{\left( x \right)}^{\varphi }},$$
where ${{\mathcal{B}}_{r}}\left( x \right)$ is the generalized binomial series [9, p.200]:
$${{\mathcal{B}}_{r}}{{\left( x \right)}^{\varphi }}=\sum\limits_{n=0}^{\infty }{\frac{\varphi }{\varphi +rn}}\left( \begin{matrix}
   \varphi +rn  \\
   n  \\
\end{matrix} \right){{x}^{n}}.$$
Thus,
$$\left[ {{x}^{n}} \right]{{g}^{\varphi }}\left( x \right)=\sum\limits_{v=0}^{\left\lfloor \frac{n}{2r+1} \right\rfloor }{b_{0}^{n-v\left( 2r+1 \right)}}b_{r}^{v}\left( \begin{matrix}
   \varphi +v\left( r+1 \right)+n-v\left( 2r+1 \right)-1  \\
   n-v\left( 2r+1 \right)  \\
\end{matrix} \right)\left[ {{x}^{v}} \right]{{\mathcal{B}}_{r+1}}{{\left( x \right)}^{\varphi }},$$
or
$${{g}^{\varphi }}\left( x \right)=\left( \frac{1}{{{\left( 1-{{b}_{0}}x \right)}^{\varphi }}},\frac{{{b}_{r}}{{x}^{2r+1}}}{{{\left( 1-{{b}_{0}}x \right)}^{r+1}}} \right){{\mathcal{B}}_{r+1}}{{\left( x \right)}^{\varphi }}.$$
In particular,
$$\left( \frac{1}{{{\left( 1-{{b}_{0}}x \right)}^{\varphi }}},\frac{{{b}_{1}}{{x}^{3}}}{{{\left( 1-{{b}_{0}}x \right)}^{2}}} \right){{\left( \frac{1-\sqrt{1-4x}}{2x} \right)}^{\varphi }}
={{\left( \frac{1-{{b}_{0}}x-\sqrt{{{\left( 1-{{b}_{0}}x \right)}^{2}}-4{{b}_{1}}{{x}^{3}}}}{2{{b}_{1}}{{x}^{3}}} \right)}^{\varphi }}.$$
{\bfseries Example 3.2.}
This example is related to the previous example and emphasizes the analogy between formulas (1) and (2). Let $A\left( x \right)=1+{{a}_{1}}x+{{a}_{r}}{{x}^{r}}$. The monomial $a_{1}^{n-vr}a_{r}^{v}$, $v=0$, $1$, … $\left\lfloor {n}/{r}\; \right\rfloor $, corresponds to the partition of  number $n$ into  parts equal to $1$ and $r$.  The coefficient of this monomial in (1) is
$$\frac{\varphi {{\left( \varphi +n-1 \right)}_{n-vr+v-1}}}{\left( n-vr \right)!v!}=\frac{\varphi \left( \varphi +n-1 \right)!}{\left( \varphi +v\left( r-1 \right) \right)!\left( n-vr \right)!v!}=$$
$$=\left( \begin{matrix}
   \varphi +n-1  \\
   n-vr  \\
\end{matrix} \right)\frac{\varphi }{\varphi +rv}\left( \begin{matrix}
   \varphi +rv  \\
   v  \\
\end{matrix} \right)=\left( \begin{matrix}
   \varphi +vr+n-vr-1  \\
   n-vr  \\
\end{matrix} \right)\left[ {{x}^{v}} \right]{{\mathcal{B}}_{r}}{{\left( x \right)}^{\varphi }}.$$
Thus,
$$\left[ {{x}^{n}} \right]{{g}^{\varphi }}\left( x \right)=\sum\limits_{v=0}^{\left\lfloor \frac{n}{r} \right\rfloor }{a_{1}^{n-vr}}a_{r}^{v}\left( \begin{matrix}
   \varphi +vr+n-vr-1  \\
   n-vr  \\
\end{matrix} \right)\left[ {{x}^{v}} \right]{{\mathcal{B}}_{r}}{{\left( x \right)}^{\varphi }},$$
$${{g}^{\varphi }}\left( x \right)=\left( \frac{1}{{{\left( 1-{{a}_{1}}x \right)}^{\varphi }}},\frac{{{a}_{r}}{{x}^{r}}}{{{\left( 1-{{a}_{1}}x \right)}^{r}}} \right){{\mathcal{B}}_{r}}{{\left( x \right)}^{\varphi }}.$$
In particular,
$$\left( \frac{1}{{{\left( 1-{{a}_{1}}x \right)}^{\varphi }}},\frac{{{a}_{2}}{{x}^{2}}}{{{\left( 1-{{a}_{1}}x \right)}^{2}}} \right){{\left( \frac{1-\sqrt{1-4x}}{2x} \right)}^{\varphi }}
={{\left( \frac{1-{{a}_{1}}x-\sqrt{{{\left( 1-{{a}_{1}}x \right)}^{2}}-4{{a}_{2}}{{x}^{2}}}}{2{{a}_{2}}{{x}^{2}}} \right)}^{\varphi }}.$$
\section{Composition polynomials of  RNA matrix}
Let $\left( R\left( x \right),xR\left( x \right) \right)$ be the RNA matrix (А097724, [20]):
$$\left( R\left( x \right),xR\left( x \right) \right)=\left( \begin{matrix}
   1 & 0 & 0 & 0 & 0 & 0 & 0 & \cdots   \\
   1 & 1 & 0 & 0 & 0 & 0 & 0 & \cdots   \\
   1 & 2 & 1 & 0 & 0 & 0 & 0 & \cdots   \\
   2 & 3 & 3 & 1 & 0 & 0 & 0 & \cdots   \\
   4 & 6 & 6 & 4 & 1 & 0 & 0 & \cdots   \\
   8 & 13 & 13 & 10 & 5 & 1 & 0 & \cdots   \\
   17 & 28 & 30 & 24 & 15 & 6 & 1 & \cdots   \\
   \vdots  & \vdots  & \vdots  & \vdots  & \vdots  & \vdots  & \vdots  & \ddots   \\
\end{matrix} \right),$$
$${{\left( R\left( x \right),xR\left( x \right) \right)}^{\varphi }}={{\left( C\left( {{x}^{2}} \right),xC\left( {{x}^{2}} \right) \right)}^{-1}}{{P}^{\varphi }}\left( C\left( {{x}^{2}} \right),xC\left( {{x}^{2}} \right) \right)=$$
$$=\left( \frac{1}{1+{{x}^{2}}},\frac{x}{1+{{x}^{2}}} \right)\left( \frac{1}{1-\varphi x},\frac{x}{1-\varphi x} \right)\left( \frac{1-\sqrt{1-4{{x}^{2}}}}{2{{x}^{2}}},\frac{1-\sqrt{1-4{{x}^{2}}}}{2x} \right),$$
$${{R}^{\left( \varphi  \right)}}\left( x \right)=\frac{1-\varphi x+{{x}^{2}}-\sqrt{{{\left( 1-\varphi x+{{x}^{2}} \right)}^{2}}-4{{x}^{2}}}}{2{{x}^{2}}}.$$
The name of this matrix is due to the fact that the $n$-th coefficient of  the series $R\left( x \right)$, $n>0$,  is equal to the number of possible secondary structures of the RNA molecule consisting of the $n$ nucleotides. The series ${{R}^{\left( \varphi  \right)}}\left( x \right)$ is the solution to the equation
$$g\left( x \right)=1+xg\left( x \right)\left( \frac{\varphi }{1-{{x}^{2}}g\left( x \right)} \right),$$
so  that  the $B$-function of the matrix ${{\left( R\left( x \right),xR\left( x \right) \right)}^{\varphi }}$ is the series $\varphi {{\left( 1-x \right)}^{-1}}$. It follows from formula (2) that the composition polynomials of the RNA matrix (we denote them ${{r}_{n}}\left( x \right)$) have the form
$${{r}_{0}}\left( x \right)=1,  \qquad{{r}_{n}}\left( x \right)=\sum\limits_{m=1}^{n}{\left( \sum\limits_{n,m}{\frac{{{\left( \frac{n+m}{2} \right)}_{m-1}}}{{{m}_{0}}!{{m}_{1}}!...{{m}_{p}}!}} \right){{x}^{m}}},$$
 where the summation of the coefficient of ${{x}^{m}}$  is over all partitions $n=\sum\nolimits_{i=0}^{p}{{{m}_{i}}\left( 2i+1 \right)}$, $\sum\nolimits_{i=0}^{p}{{{m}_{i}}=m}$. Using this formula, we will begin to build the matrix $R=L\left( R\left( x \right) \right)$, $\left[ n,\to  \right]R={{r}_{n}}\left( x \right)$:
$$R=\left( \setcounter{MaxMatrixCols}{20}\begin{matrix}
   1 & 0 & 0 & 0 & 0 & 0 & 0 & 0 & 0 & 0 & 0 & \cdots   \\
   0 & 1 & 0 & 0 & 0 & 0 & 0 & 0 & 0 & 0 & 0 & \cdots   \\
   0 & 0 & 1 & 0 & 0 & 0 & 0 & 0 & 0 & 0 & 0 & \cdots   \\
   0 & 1 & 0 & 1 & 0 & 0 & 0 & 0 & 0 & 0 & 0 & \cdots   \\
   0 & 0 & 3 & 0 & 1 & 0 & 0 & 0 & 0 & 0 & 0 & \cdots   \\
   0 & 1 & 0 & 6 & 0 & 1 & 0 & 0 & 0 & 0 & 0 & \cdots   \\
   0 & 0 & 6 & 0 & 10 & 0 & 1 & 0 & 0 & 0 & 0 & \cdots   \\
   0 & 1 & 0 & 20 & 0 & 15 & 0 & 1 & 0 & 0 & 0 & \cdots   \\
   0 & 0 & 10 & 0 & 50 & 0 & 21 & 0 & 1 & 0 & 0 & \cdots   \\
   0 & 1 & 0 & 50 & 0 & 105 & 0 & 28 & 0 & 1 & 0 & \cdots   \\
   0 & 0 & 15 & 0 & 175 & 0 & 196 & 0 & 36 & 0 & 1 & \cdots   \\
   \vdots  & \vdots  & \vdots  & \vdots  & \vdots  & \vdots  & \vdots  & \vdots  & \vdots  & \vdots  & \vdots  & \ddots   \\
\end{matrix} \right).$$
The form of this matrix leads to the assumption that $\left[ 2n,\nearrow  \right]R={{N}_{n}}\left( x \right)$, where ${{N}_{n}}\left( x \right)$ are the Narayana polynomials:
$${{N}_{0}}\left( x \right)=1,  \qquad{{N}_{n}}\left( x \right)=\frac{1}{n}\sum\limits_{m=0}^{n}{\left( \begin{matrix}
   n  \\
   m-1  \\
\end{matrix} \right)\left( \begin{matrix}
   n  \\
   m  \\
\end{matrix} \right){{x}^{m}}}.$$
Let's turn to the matrix $N$ (A090181, [20]), $\left[ n,\to  \right]N={{N}_{n}}\left( x \right)$:
$$N=\left( \begin{matrix}
   1 & 0 & 0 & 0 & 0 & 0 & 0 & \cdots   \\
   0 & 1 & 0 & 0 & 0 & 0 & 0 & \cdots   \\
   0 & 1 & 1 & 0 & 0 & 0 & 0 & \cdots   \\
   0 & 1 & 3 & 1 & 0 & 0 & 0 & \cdots   \\
   0 & 1 & 6 & 6 & 1 & 0 & 0 & \cdots   \\
   0 & 1 & 10 & 20 & 10 & 1 & 0 & \cdots   \\
   0 & 1 & 15 & 50 & 50 & 15 & 1 & \cdots   \\
   \vdots  & \vdots  & \vdots  & \vdots  & \vdots  & \vdots  & \vdots  & \ddots   \\
\end{matrix} \right).$$
{\bfseries Lemma 4.1.} \emph{$$N{{x}^{n+1}}=\frac{{{x}^{n}}{{N}_{n}}\left( x \right)}{{{\left( 1-x \right)}^{2n+1}}}, \qquad n>0.$$}
{\bfseries Proof.} The generating function of the sequence of Narayana polynomials is
$$N\left( t,x \right)=\sum\limits_{n=0}^{\infty }{{{N}_{n}}}\left( t \right){{x}^{n}}=\frac{1+x\left( 1-t \right)-\sqrt{1-2x\left( 1+t \right)+{{x}^{2}}{{\left( 1-t \right)}^{2}}}}{2x}.$$
Then
$$N\left( x,t \right)=\sum\limits_{n=0}^{\infty }{{{N}_{n}}}\left( x \right){{t}^{n}}=\frac{1+t\left( 1-x \right)-\sqrt{1-2t\left( 1+x \right)+{{t}^{2}}{{\left( 1-x \right)}^{2}}}}{2t},$$
$$N\frac{1}{1-tx}=1-t+t\sum\limits_{n=0}^{\infty }{\frac{{{N}_{n}}\left( x \right){{x}^{n}}{{t}^{n}}}{{{\left( 1-x \right)}^{2n+1}}}}=1-t+\frac{t}{1-x}N\left( x,\frac{xt}{{{\left( 1-x \right)}^{2}}} \right)=$$
$$=\frac{1+x\left( 1-t \right)-\sqrt{1-2x\left( 1+t \right)+{{x}^{2}}{{\left( 1-t \right)}^{2}}}}{2x}=\sum\limits_{n=0}^{\infty }{{{N}_{n}}}\left( t \right){{x}^{n}}.\qquad\square $$
{\bfseries Theorem 4.2.} \emph{$$\left[ 2n,\nearrow  \right]R={{N}_{n}}\left( x \right).$$}
{\bfseries Proof.} Denote ${{\tilde{N}}_{0}}\left( x \right)=1$, ${{\tilde{N}}_{n}}\left( x \right)=\left( {1}/{x}\; \right){{N}_{n}}\left( x \right)$. Then
$$\tilde{N}\left( x,t \right)=\sum\limits_{n=0}^{\infty }{{{{\tilde{N}}}_{n}}}\left( x \right){{t}^{n}}=\frac{1-t\left( 1-x \right)-\sqrt{1-2t\left( 1+x \right)+{{t}^{2}}{{\left( 1-x \right)}^{2}}}}{2xt}.$$
By Lemma 4.1., if $\left[ 2n,\nearrow  \right]R={{N}_{n}}\left( x \right)$, then
$$R{{x}^{n+1}}=\frac{{{x}^{n+1}}{{{\tilde{N}}}_{n}}\left( {{x}^{2}} \right)}{{{\left( 1-{{x}^{2}} \right)}^{2n+1}}}.$$
Then
$$R\frac{1}{1-tx}=1+xt\sum\limits_{n=0}^{\infty }{\frac{{{{\tilde{N}}}_{n}}\left( {{x}^{2}} \right){{x}^{n}}{{t}^{n}}}{{{\left( 1-{{x}^{2}} \right)}^{2n+1}}}}=1+\frac{xt}{1-{{x}^{2}}}\tilde{N}\left( {{x}^{2}},\frac{xt}{{{\left( 1-{{x}^{2}} \right)}^{2}}} \right)=$$
$$=\frac{1-tx+{{x}^{2}}-\sqrt{{{\left( 1-tx+{{x}^{2}} \right)}^{2}}-4{{x}^{2}}}}{2{{x}^{2}}}={{R}^{\left( t \right)}}\left( x \right).\qquad\square $$
Thus,
$${{r}_{2n}}\left( x \right)=\sum\limits_{m=0}^{n}{{{N}_{n+m,2m}}}{{x}^{2m}},  \qquad{{r}_{2n+1}}\left( x \right)=\sum\limits_{m=0}^{n}{{{N}_{n+m+1,2m+1}}}{{x}^{2m+1}},$$
where ${{N}_{n}}\left( x \right)=\sum\nolimits_{m=0}^{n}{{{N}_{n,m}}}{{x}^{m}}$, or 
$${{r}_{2n}}\left( x \right)=\sum\limits_{m=1}^{n}{\frac{1}{n+m}\left( \begin{matrix}
   n+m  \\
   2m-1  \\
\end{matrix} \right)\left( \begin{matrix}
   n+m  \\
   2m  \\
\end{matrix} \right)}{{x}^{2m}},$$
$${{r}_{2n+1}}\left( x \right)=\sum\limits_{m=0}^{n}{\frac{1}{n+m+1}\left( \begin{matrix}
   n+m+1  \\
   2m  \\
\end{matrix} \right)\left( \begin{matrix}
   n+m+1  \\
   2m+1  \\
\end{matrix} \right)}{{x}^{2m+1}}.$$

A generalization of the RNA matrix is the matrix $\left( R\left( \beta ,x \right),xR\left( \beta ,x \right) \right)$:
$${{\left( R\left( \beta ,x \right),xR\left( \beta ,x \right) \right)}^{\varphi }}={{\left( C\left( \beta {{x}^{2}} \right),xC\left( \beta {{x}^{2}} \right) \right)}^{-1}}{{P}^{\varphi }}\left( C\left( \beta {{x}^{2}} \right),xC\left( \beta {{x}^{2}} \right) \right)=$$
$$=\left( \frac{1}{1+\beta {{x}^{2}}},\frac{x}{1+\beta {{x}^{2}}} \right)\left( \frac{1}{1-\varphi x},\frac{x}{1-\varphi x} \right)\left( \frac{1-\sqrt{1-4\beta {{x}^{2}}}}{2\beta {{x}^{2}}},\frac{1-\sqrt{1-4\beta {{x}^{2}}}}{2\beta x} \right),$$
$${{R}^{\left( \varphi  \right)}}\left( \beta ,x \right)=\frac{1-\varphi x+\beta {{x}^{2}}-\sqrt{{{\left( 1-\varphi x+\beta {{x}^{2}} \right)}^{2}}-4\beta {{x}^{2}}}}{2\beta {{x}^{2}}}.$$
The series ${{R}^{\left( \varphi  \right)}}\left( \beta ,x \right)$ is the solution to the equation
$$g\left( x \right)=1+xg\left( x \right)\left( \frac{\varphi }{1-\beta {{x}^{2}}g\left( x \right)} \right),$$
so  that the $B$-function of the matrix ${{\left( R\left( \beta ,x \right),xR\left( \beta ,x \right) \right)}^{\varphi }}$  is the series $\varphi {{\left( 1-\beta x \right)}^{-1}}$. Hence, the composition polynomials of the matrix $\left( R\left( \beta ,x \right),xR\left( \beta ,x \right) \right)$ have the form
$$\sum\limits_{m=1}^{n}{\left( \sum\limits_{n,m}{\frac{{{\left( \frac{n+m}{2} \right)}_{m-1}}}{{{m}_{0}}!{{m}_{1}}!...{{m}_{p}}!}} \right){{\beta }^{\frac{n-m}{2}}}{{x}^{m}}={{\left( \sqrt{\beta } \right)}^{n}}{{r}_{n}}\left( {x}/{\sqrt{\beta }}\; \right)}.$$
\section{$B$-composition  polynomials}
The Bell subgroup matrix whose $B$-sequence has the generating function $\varphi B\left( x \right)$ will be denoted by ${{\left( g\left( x \right),xg\left( x \right) \right)}^{\left[ \varphi  \right]}}=\left( {{g}^{\left[ \varphi  \right]}}\left( x \right),x{{g}^{\left[ \varphi  \right]}}\left( x \right) \right)$. The polynomials  ${{u}_{n}}\left( x \right)$ such that ${{g}^{\left[ \varphi  \right]}}\left( x \right)=\sum\nolimits_{n=0}^{\infty }{{{u}_{n}}}\left( \varphi  \right){{x}^{n}}$  will be called $B$-composition polynomials. It follows from formula (2) that 
$$\left[ {{x}^{m}} \right]{{u}_{n}}\left( x \right)={{\left( \frac{n+m}{2} \right)}_{m-1}}\sum\limits_{n,m}{\frac{b_{0}^{{{m}_{0}}}b_{1}^{{{m}_{1}}}...b_{p}^{{{m}_{p}}}}{{{m}_{0}}!{{m}_{1}}!...{{m}_{p}}!}},\eqno (3)$$
where the summation is over all partitions $n=\sum\nolimits_{i=0}^{p}{{{m}_{i}}\left( 2i+1 \right)}$, $\sum\nolimits_{i=0}^{p}{{{m}_{i}}=m}$. It is also seen from formula (2) that if the $B$-function $B\left( x \right)$ is related to the polynomials ${{u}_{n}}\left( x \right)$, then the $B$-function $B\left( \beta x \right)$ is related to the polynomials ${{\left( \sqrt{\beta } \right)}^{n}}{{u}_{n}}\left( {x}/{\sqrt{\beta }}\; \right)$.

The $B$-expansion when $v=1$ will be called ${{B}_{1}}$-expansion. The initial terms of the ${{B}_{1}}$-expansion are:
$${{g}_{0}}=1,  \quad{{g}_{1}}={{b}_{0}},  \quad{{g}_{2}}=b_{0}^{2},  \quad{{g}_{3}}=b_{0}^{3}+{{b}_{1}},  \quad{{g}_{4}}=b_{0}^{4}+3{{b}_{0}}{{b}_{1}}, $$
$${{g}_{5}}=b_{0}^{5}+6b_{0}^{2}{{b}_{1}}+{{b}_{2}},  \qquad{{g}_{6}}=b_{0}^{5}+10b_{0}^{3}{{b}_{1}}+4{{b}_{0}}{{b}_{2}}+2b_{1}^{2},$$ 
$${{g}_{7}}=b_{0}^{7}+15b_{0}^{4}{{b}_{1}}+10b_{0}^{2}{{b}_{2}}+10{{b}_{0}}b_{1}^{2}+{{b}_{3}},$$
$${{g}_{8}}=b_{0}^{8}+21b_{0}^{5}{{b}_{1}}+20b_{0}^{3}{{b}_{2}}+30b_{0}^{2}b_{1}^{2}+5{{b}_{0}}{{b}_{3}}+5{{b}_{1}}{{b}_{2}},$$
$${{g}_{9}}=b_{0}^{9}+28b_{0}^{6}{{b}_{1}}+35b_{0}^{4}{{b}_{2}}+70b_{0}^{3}b_{1}^{2}+15b_{0}^{2}{{b}_{3}}+30{{b}_{0}}{{b}_{1}}{{b}_{2}}+5b_{1}^{3}+{{b}_{4}},$$
$${{g}_{10}}=b_{0}^{10}+36b_{0}^{7}{{b}_{1}}+56b_{0}^{5}{{b}_{2}}+140b_{0}^{4}b_{1}^{2}+35b_{0}^{3}{{b}_{3}}+35{{b}_{0}}b_{1}^{3}+105b_{0}^{2}{{b}_{1}}{{b}_{2}}+$$
$$+6{{b}_{0}}{{b}_{4}}+6{{b}_{1}}{{b}_{3}}+3b_{2}^{2}.$$

Using formula (3), we can build the matrix whose rows correspond to the $B$-compositions polynomials. We call such matrix a $B$-composition matrix. Note that the first column (i.e., $n=1$) of such matrix has the generating function $xB\left( {{x}^{2}} \right)$. 
\section{Case $B\left( x \right)=1+x$}
The series 
$$_{\left( 1 \right)}{{R}^{\left[ \varphi  \right]}}\left( x \right)=\frac{1-\varphi x-\sqrt{{{\left( 1-\varphi x \right)}^{2}}-4\varphi {{x}^{3}}}}{2\varphi {{x}^{3}}}$$
is the solution to the equation $g\left( x \right)=1+xg\left( x \right)\varphi \left( 1+{{x}^{2}}g\left( x \right) \right)$, so  that the $B$-function of the matrix ${{\left( {}_{\left( 1 \right)}R\left( x \right),x{}_{\left( 1 \right)}R\left( x \right) \right)}^{\left[ \varphi  \right]}}$ is the series $\varphi \left( 1+x \right)$. The $B$-composition matrix has the form
$$_{\left( 1 \right)}R=\left( \begin{matrix}
   1 & 0 & 0 & 0 & 0 & 0 & 0 & 0 & 0 & 0 & 0 & \cdots   \\
   0 & 1 & 0 & 0 & 0 & 0 & 0 & 0 & 0 & 0 & 0 & \cdots   \\
   0 & 0 & 1 & 0 & 0 & 0 & 0 & 0 & 0 & 0 & 0 & \cdots   \\
   0 & 1 & 0 & 1 & 0 & 0 & 0 & 0 & 0 & 0 & 0 & \cdots   \\
   0 & 0 & 3 & 0 & 1 & 0 & 0 & 0 & 0 & 0 & 0 & \cdots   \\
   0 & 0 & 0 & 6 & 0 & 1 & 0 & 0 & 0 & 0 & 0 & \cdots   \\
   0 & 0 & 2 & 0 & 10 & 0 & 1 & 0 & 0 & 0 & 0 & \cdots   \\
   0 & 0 & 0 & 10 & 0 & 15 & 0 & 1 & 0 & 0 & 0 & \cdots   \\
   0 & 0 & 0 & 0 & 30 & 0 & 21 & 0 & 1 & 0 & 0 & \cdots   \\
   0 & 0 & 0 & 5 & 0 & 70 & 0 & 28 & 0 & 1 & 0 & \cdots   \\
   0 & 0 & 0 & 0 & 35 & 0 & 140 & 0 & 36 & 0 & 1 & \cdots   \\
   \vdots  & \vdots  & \vdots  & \vdots  & \vdots  & \vdots  & \vdots  & \vdots  & \vdots  & \vdots  & \vdots  & \ddots   \\
\end{matrix} \right).$$
The coefficient of monomial $b_{0}^{p}b_{1}^{v}$ in the ${{B}_{1}}$-expansion is 
$$\frac{{{\left( p+2v \right)}_{p+v-1}}}{p!v!}=\frac{\left( p+2v \right)!}{\left( 1+v \right)!p!v!}=\frac{1}{1+v}\left( \begin{matrix}
   p+2v  \\
   p  \\
\end{matrix} \right)\left( \begin{matrix}
   2v  \\
   v  \\
\end{matrix} \right)={{C}_{v}}\left( \begin{matrix}
   2v+p  \\
   p  \\
\end{matrix} \right),$$ 
${{C}_{v}}=\left[ {{x}^{v}} \right]C\left( x \right)$. The monomial $b_{0}^{p}b_{1}^{v}$ corresponds to the partition of number $n=p+3v$ into $m=p+v$ parts. Hence,
$${{\left( _{\left( 1 \right)}R \right)}_{n,m}}={{C}_{{\left( n-m \right)}/{2}\;}}\left( \begin{matrix}
   {\left( n+m \right)}/{2}\;  \\
   {\left( 3m-n \right)}/{2}\;  \\
\end{matrix} \right),$$
where ${{C}_{{\left( n-m \right)}/{2}\;}}=0$, if $n-m$ is odd,
$$\left[ 2n,\searrow  \right]{}_{\left( 1 \right)}R=\sum\limits_{m=0}^{\infty }{{{C}_{n}}}\left( \begin{matrix}
   n+m  \\
   m-n  \\
\end{matrix} \right){{x}^{m}}={{x}^{n}}{{C}_{n}}{{\left( \frac{1}{1-x} \right)}^{2n+1}},$$
$$\left[ 2n,\nearrow  \right]{}_{\left( 1 \right)}R=\sum\limits_{m=0}^{n}{{{C}_{n-m}}}\left( \begin{matrix}
   n  \\
   2m-n  \\
\end{matrix} \right){{x}^{m}},$$
$${}_{\left( 1 \right)}R{{x}^{2n}}={{x}^{2n}}\sum\limits_{m=0}^{n}{{{C}_{m}}\left( \begin{matrix}
   2n+m  \\
   2n-m  \\
\end{matrix} \right)}{{x}^{2m}},  \quad{}_{\left( 1 \right)}R{{x}^{2n+1}}={{x}^{2n+1}}\sum\limits_{m=0}^{n}{{{C}_{m}}\left( \begin{matrix}
   2n+1+m  \\
   2n+1-m  \\
\end{matrix} \right)}{{x}^{2m}},$$
$$\left[ 2n,\to  \right]{}_{\left( 1 \right)}R={}_{\left( 1 \right)}{{r}_{2n}}\left( x \right)=\sum\limits_{m=0}^{n}{{{C}_{n-m}}\left( \begin{matrix}
   n+m  \\
   3m-n  \\
\end{matrix} \right)}{{x}^{2m}},$$
$$\left[ 2n+1,\to  \right]{}_{\left( 1 \right)}R={}_{\left( 1 \right)}{{r}_{2n+1}}\left( x \right)=\sum\limits_{m=0}^{n}{{{C}_{n-m}}\left( \begin{matrix}
   n+1+m  \\
   3m+1-n  \\
\end{matrix} \right)}{{x}^{2m+1}}.$$
Let's turn to the polynomials ${{T}_{n}}\left( x \right)$ ( A033282, [20]):
$${{T}_{n}}\left( x \right)=\frac{1}{n+1}\sum\limits_{m=0}^{n}{\left( \begin{matrix}
   n+1  \\
   m+1  \\
\end{matrix} \right)\left( \begin{matrix}
   n+m+2  \\
   m  \\
\end{matrix} \right){{x}^{m}}}={{\left( 1+x \right)}^{n}}{{\tilde{N}}_{n+1}}\left( \frac{x}{1+x} \right).$$
Since
$$\sum\limits_{n=0}^{\infty }{{{{\tilde{N}}}_{n+1}}}\left( x \right){{t}^{n}}=\frac{1-t\left( 1+x \right)-\sqrt{1-2t\left( 1+x \right)+{{t}^{2}}{{\left( 1-x \right)}^{2}}}}{2x{{t}^{2}}}=\bar{N}\left( x,t \right),$$
then
$$T\left( x,t \right)=\sum\limits_{n=0}^{\infty }{{{T}_{n}}}\left( x \right){{t}^{n}}=\bar{N}\left( \frac{x}{1+x},\left( 1+x \right)t \right)=
\frac{1-t\left( 1+2x \right)-\sqrt{1-2t\left( 1+2x \right)+{{t}^{2}}}}{2x{{t}^{2}}\left( 1+x \right)}.$$
\\{\bfseries Theorem 6.1.} \emph{$$_{\left( 1 \right)}R{{x}^{n+1}}={{x}^{n+1}}{{T}_{n}}\left( {{x}^{2}} \right)\left( 1+{{x}^{2}} \right).$$}
{\bfseries Proof.}  
$$_{\left( 1 \right)}R\frac{1}{1-tx}=1+tx\left( 1+{{x}^{2}} \right)\sum\limits_{n=0}^{\infty }{{{T}_{n}}}\left( {{x}^{2}} \right){{x}^{n}}{{t}^{n}}=1+tx\left( 1+{{x}^{2}} \right)T\left( {{x}^{2}},xt \right)=$$
$$=\frac{1-tx-\sqrt{{{\left( 1-tx \right)}^{2}}-4t{{x}^{3}}}}{2t{{x}^{3}}}={}_{\left( 1 \right)}{{R}^{\left[ t \right]}}\left( x \right).\qquad \square$$
\section{Case $B\left( x \right)=C\left( x \right)$}
The series 
$$_{\left( 2 \right)}{{R}^{\left[ \varphi  \right]}}\left( x \right)=\frac{1+\left( \left( {2}/{\varphi }\; \right)-\varphi  \right)x-\sqrt{1-2\varphi x+\left( {{\varphi }^{2}}-4 \right){{x}^{2}}}}{2x\left( {1}/{\varphi }\; \right)}$$
is the solution to the equation 
$$g\left( x \right)=1+xg\left( x \right)\varphi \left( \frac{1-\sqrt{1-4{{x}^{2}}g\left( x \right)}}{2{{x}^{2}}g\left( x \right)} \right),$$
so  that the $B$-function of the matrix ${{\left( {}_{\left( 2 \right)}R\left( x \right),x{}_{\left( 2 \right)}R\left( x \right) \right)}^{\left[ \varphi  \right]}}$  is the series $\varphi C\left( x \right)$. The $B$-composition matrix has the form
$$_{\left( 2 \right)}R=\left( \begin{matrix}
   1 & 0 & 0 & 0 & 0 & 0 & 0 & 0 & 0 & 0 & 0 & \cdots   \\
   0 & 1 & 0 & 0 & 0 & 0 & 0 & 0 & 0 & 0 & 0 & \cdots   \\
   0 & 0 & 1 & 0 & 0 & 0 & 0 & 0 & 0 & 0 & 0 & \cdots   \\
   0 & 1 & 0 & 1 & 0 & 0 & 0 & 0 & 0 & 0 & 0 & \cdots   \\
   0 & 0 & 3 & 0 & 1 & 0 & 0 & 0 & 0 & 0 & 0 & \cdots   \\
   0 & 2 & 0 & 6 & 0 & 1 & 0 & 0 & 0 & 0 & 0 & \cdots   \\
   0 & 0 & 10 & 0 & 10 & 0 & 1 & 0 & 0 & 0 & 0 & \cdots   \\
   0 & 5 & 0 & 30 & 0 & 15 & 0 & 1 & 0 & 0 & 0 & \cdots   \\
   0 & 0 & 35 & 0 & 70 & 0 & 21 & 0 & 1 & 0 & 0 & \cdots   \\
   0 & 14 & 0 & 140 & 0 & 140 & 0 & 28 & 0 & 1 & 0 & \cdots   \\
   0 & 0 & 126 & 0 & 420 & 0 & 252 & 0 & 36 & 0 & 1 & \cdots   \\
   \vdots  & \vdots  & \vdots  & \vdots  & \vdots  & \vdots  & \vdots  & \vdots  & \vdots  & \vdots  & \vdots  & \ddots   \\
\end{matrix} \right).$$
We assume that $\left[ 2n,\searrow  \right]{}_{\left( 2 \right)}R=\left( {1}/{{{x}^{n-1}}}\; \right)\left[ 2n,\searrow  \right]{}_{\left( 1 \right)}R$, $n>0$. Let's turn to the matrix $_{\left( 1,2 \right)}R$ (А107131, [20]),  $_{\left( 1,2 \right)}R{{x}^{n+1}}={{x}^{n+1}}{{T}_{n}}\left( x \right)\left( 1+x \right)$:
$$_{\left( 1,2 \right)}R=\left( \begin{matrix}
   1 & 0 & 0 & 0 & 0 & 0 & 0 & 0 & \cdots   \\
   0 & 1 & 0 & 0 & 0 & 0 & 0 & 0 & \cdots   \\
   0 & 1 & 1 & 0 & 0 & 0 & 0 & 0 & \cdots   \\
   0 & 0 & 3 & 1 & 0 & 0 & 0 & 0 & \cdots   \\
   0 & 0 & 2 & 6 & 1 & 0 & 0 & 0 & \cdots   \\
   0 & 0 & 0 & 10 & 10 & 1 & 0 & 0 & \cdots   \\
   0 & 0 & 0 & 5 & 30 & 15 & 1 & 0 & \cdots   \\
   0 & 0 & 0 & 0 & 35 & 70 & 21 & 1 & \cdots   \\
   \vdots  & \vdots  & \vdots  & \vdots  & \vdots  & \vdots  & \vdots  & \vdots  & \ddots   \\
\end{matrix} \right).$$
Denote $\left[ n,\to  \right]{}_{\left( 1,2 \right)}R={{F}_{n}}\left( x \right)$. Then
$${}_{\left( 1,2 \right)}R\frac{1}{1-tx}=F\left( t,x \right)=\sum\limits_{n=0}^{\infty }{{{F}_{n}}}\left( t \right){{x}^{n}}=1+xt\left( 1+x \right)\sum\limits_{n=0}^{\infty }{{{T}_{n}}}\left( x \right){{x}^{n}}{{t}^{n}}=$$
$$=\frac{1-xt-\sqrt{1-2xt\left( 1+2x \right)+{{x}^{2}}{{t}^{2}}}}{2{{x}^{2}}t}.$$
{\bfseries Theorem 7.1.} \emph{$$\left[ n+1,\to  \right]{}_{\left( 2 \right)}R=\frac{1}{{{x}^{n-1}}}{{F}_{n}}\left( {{x}^{2}} \right).$$}
{\bfseries Proof.} 
$${}_{\left( 2 \right)}R\frac{1}{1-tx}=1+xt\sum\limits_{n=0}^{\infty }{{{F}_{n}}}\left( {{t}^{2}} \right)\frac{{{x}^{n}}}{{{t}^{n}}}=1+xtF\left( {{t}^{2}},{x}/{t}\; \right)=$$
$$=\frac{1+\left( \left( {2}/{t}\; \right)-t \right)x-\sqrt{1-2tx+\left( {{t}^{2}}-4 \right){{x}^{2}}}}{2x\left( {1}/{t}\; \right)}={}_{\left( 2 \right)}{{R}^{\left[ t \right]}}\left( x \right).\qquad \square$$
Thus,
$${{\left( _{\left( 2 \right)}R \right)}_{n,m}}={{C}_{{\left( n-m \right)}/{2}\;}}\left( \begin{matrix}
   n-1  \\
   m-1  \\
\end{matrix} \right),$$
$$\left[ 2n,\to  \right]{}_{\left( 2 \right)}R={}_{\left( 2 \right)}{{r}_{2n}}\left( x \right)=\sum\limits_{m=0}^{n}{{{C}_{n-m}}\left( \begin{matrix}
   2n-1  \\
   2m-1  \\
\end{matrix} \right)}{{x}^{2m}},$$
$$\left[ 2n+1,\to  \right]{}_{\left( 2 \right)}R={}_{\left( 2 \right)}{{r}_{2n+1}}\left( x \right)=\sum\limits_{m=0}^{n}{{{C}_{n-m}}\left( \begin{matrix}
   2n  \\
   2m  \\
\end{matrix} \right)}{{x}^{2m+1}}.$$

Denote $\left( {1}/{x}\; \right){}_{\left( 2 \right)}{{r}_{n+1}}\left( x \right)={}_{\left( 2 \right)}{{\tilde{r}}_{n}}\left( x \right)$. Since 
$${{\left( x,x \right)}^{T}}{}_{\left( 2 \right)}R\left( x,x \right)={{\left( \tilde{C}\left( x \right),x \right)}_{E}}, \qquad\tilde{C}\left( x \right)=\sum\limits_{n=0}^{\infty }{{{C}_{n}}}\frac{{{x}^{2n}}}{\left( 2n \right)!},$$
then the sequence of polynomials ${}_{\left( 2 \right)}{{\tilde{r}}_{n}}\left( x \right)$ is the Appell sequence:
$$\sum\limits_{n=0}^{\infty }{\frac{{}_{\left( 2 \right)}{{{\tilde{r}}}_{n}}\left( \varphi  \right)}{n!}}{{x}^{n}}=\tilde{C}\left( x \right){{e}^{\varphi x}}.$$
Thus,
$$\left[ {{x}^{n}} \right]{}_{\left( 2 \right)}{{R}^{\left[ \varphi  \right]}}\left( x \right)=\varphi \left( n-1 \right)!\left[ {{x}^{n-1}} \right]\tilde{C}\left( x \right){{e}^{\varphi x}}.$$

The $B$-composition matrix will be denoted by $\left\langle B\left( x \right) \right\rangle $. If
$${{\left( x,x \right)}^{T}}\left\langle B\left( x \right) \right\rangle \left( x,x \right)={{\left( \tilde{B}\left( x \right),x \right)}_{E}}, \qquad\tilde{B}\left( x \right)=\sum\limits_{n=0}^{\infty }{{{b}_{n}}}\frac{{{x}^{2n}}}{\left( 2n \right)!},$$
the matrix $\left\langle B\left( x \right) \right\rangle $ we call Appell type matrix.\\
{\bfseries Theorem 7.2.} \emph{If the matrix $\left\langle B\left( x \right) \right\rangle $ is the Appell type matrixя, then ${{b}_{n}}={{C}_{n}}b_{1}^{n}$.}\\
{\bfseries Proof.} If the matrix $\left\langle B\left( x \right) \right\rangle $ is the Appell type matrix, then the following identity holds
$$\sum\limits_{m=0}^{n}{{{b}_{n-m}}\left( \begin{matrix}
   2n+1  \\
   2m+1  \\
\end{matrix} \right)}=\sum\limits_{2\left( n+1 \right)}^{{}}{\frac{{{\left( k \right)}_{q-1}}}{{{m}_{0}}!{{m}_{1}}!...{{m}_{n}}!}}b_{0}^{{{m}_{0}}}b_{1}^{{{m}_{1}}}...b_{n}^{{{m}_{n}}},$$
  $$k=\sum\limits_{i=0}^{n}{\left( i+1 \right){{m}_{i}}}, \qquad q=\sum\limits_{i=0}^{n}{{{m}_{i}}},$$
where in the right part the summation is over all monomials $b_{0}^{{{m}_{0}}}b_{1}^{{{m}_{1}}}...b_{n}^{{{m}_{n}}}$  for which $2\left( n+1 \right)=\sum\nolimits_{i=0}^{n}{{{m}_{i}}\left( 2i+1 \right)}$. We will consider this identity as the equation with unknowns ${{b}_{1}}$, ${{b}_{2}}$ … ${{b}_{n}}$ (it's obvious that ${{b}_{0}}=1$). Since the monomial ${{b}_{0}}{{b}_{n}}$ corresponds to the partition of number $2\left( n+1 \right)$ into two parts equal to $1$ and $2n+1$, the equation can be represented as
$$\left( 2n+1 \right){{b}_{n}}+{{f}_{1}}\left( {{b}_{1}},{{b}_{2}},...,{{b}_{n-1}} \right)=\left( n+2 \right){{b}_{n}}+{{f}_{2}}\left( {{b}_{1}},{{b}_{2}},...,{{b}_{n-1}} \right),   \qquad n>0,$$
where ${{f}_{1}}\left( {{b}_{1}},{{b}_{2}},...,{{b}_{n-1}} \right)$, ${{f}_{2}}\left( {{b}_{1}},{{b}_{2}},...,{{b}_{n-1}} \right)$ are independent of ${{b}_{n}}$. Thus, the $n$th term of the $B$-sequence, starting from the second, is uniquely expressed through the previous terms. The result is known: ${{b}_{n}}={{C}_{n}}b_{1}^{n}$. \qquad  $\square $
\section{Connection theorem }
The generating function of the $n$th descending diagonal of the exponential Riordan matrix ${{\left( f\left( x \right),xg\left( x \right) \right)}_{E}}$, ${{g}_{0}}=1$, has the form ${{{h}_{n}}\left( x \right)}/{{{\left( 1-x \right)}^{2n+1}}}\;$, where ${{h}_{n}}\left( x \right)$ is a polynomial of degree $\le n$ [4,15]. In particular,
$$\left[ n,\searrow  \right]{{\left( 1,\frac{x}{1-x} \right)}_{E}}=\frac{\left( n+1 \right)!{{N}_{n}}\left( x \right)}{{{\left( 1-x \right)}^{2n+1}}},$$
$$\left[ n,\searrow  \right]{{\left( 1,x\left( 1+x \right) \right)}_{E}}=\frac{\left( {\left( 2n \right)!}/{n!}\; \right){{x}^{n}}}{{{\left( 1-x \right)}^{2n+1}}},  \quad\left[ n,\searrow  \right]{{\left( 1,xC\left( x \right) \right)}_{E}}=\frac{\left( {\left( 2n \right)!}/{n!}\; \right)x}{{{\left( 1-x \right)}^{2n+1}}},$$
(in the latter case $n>0$). Thus,
$$\left[ 2n,\searrow  \right]\left\langle \frac{1}{1-x} \right\rangle =\frac{1}{\left( n+1 \right)!}\left[ n,\searrow  \right]{{\left( 1,\frac{x}{1-x} \right)}_{E}},$$
$$\left[ 2n,\searrow  \right]\left\langle 1+x \right\rangle =\frac{1}{\left( n+1 \right)!}\left[ n,\searrow  \right]{{\left( 1,x\left( 1+x \right) \right)}_{E}},$$
$$\left[ 2n,\searrow  \right]\left\langle C\left( x \right) \right\rangle =\frac{1}{\left( n+1 \right)!}\left[ n,\searrow  \right]{{\left( 1,xC\left( x \right) \right)}_{E}}.$$
This observation leads to the following theorem.
\\{\bfseries Theorem 8.1.} \emph{$$\left[ 2n,\searrow  \right]\left\langle B\left( x \right) \right\rangle =\frac{1}{\left( n+1 \right)!}\left[ n,\searrow  \right]{{\left( 1,xB\left( x \right) \right)}_{E}}.$$}
{\bfseries Proof.}
Elements of the  matrix ${{\left( 1,xB\left( x \right) \right)}_{E}}$ are expressed in terms of  coefficients of the series $B\left( x \right)$ by the formula
$${{\left( {{\left( 1,xB\left( x \right) \right)}_{E}} \right)}_{n,m}}=n!\sum\limits_{n,m}{\frac{b_{0}^{{{m}_{0}}}b_{1}^{{{m}_{1}}}...b_{n-1}^{{{m}_{n-1}}}}{{{m}_{0}}!{{m}_{1}}!...{{m}_{n-1}}!}},$$
where the summation is over all monomials $b_{0}^{{{m}_{0}}}b_{1}^{{{m}_{1}}}...b_{n-1}^{{{m}_{n-1}}}$  for which $n=\sum\nolimits_{i=0}^{n-1}{{{m}_{i}}\left( i+1 \right)}$, $m=\sum\nolimits_{i=0}^{n-1}{{{m}_{i}}}$. A similar formula for the matrix $\left\langle B\left( x \right) \right\rangle $ has the form
$${{\left( \left\langle B\left( x \right) \right\rangle  \right)}_{n,m}}={{\left( \frac{n+m}{2} \right)}_{m-1}}\sum\limits_{n,m}{\frac{b_{0}^{{{m}_{0}}}b_{1}^{{{m}_{1}}}...b_{p}^{{{m}_{p}}}}{{{m}_{0}}!{{m}_{1}}!...{{m}_{p}}!}},$$
where the summation is over all monomials $b_{0}^{{{m}_{0}}}b_{1}^{{{m}_{1}}}...b_{p}^{{{m}_{p}}}$  for which $n=\sum\nolimits_{i=0}^{p}{{{m}_{i}}\left( 2i+1 \right)}$,  $m=\sum\nolimits_{i=0}^{p}{{{m}_{i}}}$. We must prove that
$${{\left( \left\langle B\left( x \right) \right\rangle  \right)}_{2n-m,m}}=\frac{1}{\left( n-m+1 \right)!}{{\left( {{\left( 1,xB\left( x \right) \right)}_{E}} \right)}_{n,m}}.$$
This comes down to the proof that
$$\sum\limits_{n,m}{\frac{b_{0}^{{{m}_{0}}}b_{1}^{{{m}_{1}}}...b_{n-1}^{{{m}_{n-1}}}}{{{m}_{0}}!{{m}_{1}}!...{{m}_{n-1}}!}}=\sum\limits_{2n-m,m}{\frac{b_{0}^{{{m}_{0}}}b_{1}^{{{m}_{1}}}...b_{p}^{{{m}_{p}}}}{{{m}_{0}}!{{m}_{1}}!...{{m}_{p}}!}},$$
where on the left the summation is carried by the rule $n=\sum\nolimits_{i=0}^{n-1}{{{m}_{i}}\left( i+1 \right)}$, $m=\sum\nolimits_{i=0}^{n-1}{{{m}_{i}}}$, on the right – by the rule $2n-m=\sum\nolimits_{i=0}^{p}{{{m}_{i}}\left( 2i+1 \right)}$, $m=\sum\nolimits_{i=0}^{p}{{{m}_{i}}}$. The isomorphism between the set of  partitions of  number $n$ into $m$ parts and the set of  partitions of number $2n-m$ into $m$ odd parts (each partition $n=\sum\nolimits_{i=0}^{n-m}{{{m}_{i}}\left( i+1 \right)}$ corresponds to the partition $2n-m=\sum\nolimits_{i=0}^{n-m}{{{m}_{i}}\left( 2i+1 \right)}$, and vice versa) is the proof.  \qquad     $\square $

\section{$B$-composition-convolution polynomials }
Let ${{s}_{n}}\left( x \right)$ be the convolution polynomials of the series $B\left( x \right)$: ${{B}^{m}}\left( x \right)=\sum\nolimits_{n=0}^{\infty }{{{s}_{n}}}\left( m \right){{x}^{n}}$. Then
$${{\left( {{\left( 1,xB\left( x \right) \right)}_{E}} \right)}_{n,m}}=\frac{n!{{s}_{n-m}}\left( m \right)}{m!},$$ 
$$\left[ n,\to  \right]\left\langle B\left( x \right) \right\rangle ={{u}_{n}}\left( x \right)=\sum\limits_{m=1}^{n}{\frac{{{\left( \frac{n+m}{2} \right)}_{m-1}}{{s}_{\frac{n-m}{2}}}\left( m \right)}{m!}}{{x}^{m}},$$
where ${{s}_{{\left( n-m \right)}/{2}\;}}\left( m \right)=0$, if $n-m$ is odd,
$${{u}_{2n}}\left( x \right)=\sum\limits_{m=0}^{n}{\left( \begin{matrix}
   n+m  \\
   2m  \\
\end{matrix} \right)\frac{{{s}_{n-m}}\left( 2m \right)}{n-m+1}}{{x}^{2m}},$$
$${{u}_{2n+1}}\left( x \right)=\sum\limits_{m=0}^{n}{\left( \begin{matrix}
   n+m+1  \\
   2m+1  \\
\end{matrix} \right)\frac{{{s}_{n-m}}\left( 2m+1 \right)}{n-m+1}}{{x}^{2m+1}}.$$
{\bfseries Example 9.1.}
$$B\left( x \right)={{e}^{x}},  \qquad{{u}_{2n}}\left( x \right)=\sum\limits_{m=0}^{n}{\left( \begin{matrix}
   n+m  \\
   2m  \\
\end{matrix} \right)\frac{{{\left( 2m \right)}^{n-m}}}{\left( n-m+1 \right)!}}{{x}^{2m}},$$
$${{u}_{2n+1}}\left( x \right)=\sum\limits_{m=0}^{n}{\left( \begin{matrix}
   n+m+1  \\
   2m+1  \\
\end{matrix} \right)\frac{{{\left( 2m+1 \right)}^{n-m}}}{\left( n-m+1 \right)!}}{{x}^{2m+1}}.$$
From formula (2) and Theorem 8.1 it follows that if ${{s}_{n}}\left( x \right)$ are the convolution polynomials of  the $B$-function of the matrix $\left( g\left( x \right),xg\left( x \right) \right)$, ${{q}_{n}}\left( x \right)$ are the  convolution polynomials of the series $g\left( x \right)$, then
$${{q}_{n}}\left( x \right)={{\sum\limits_{m=1}^{n}{x\left( x+\frac{n+m}{2}-1 \right)}}_{m-1}}\frac{{{s}_{{\left( n-m \right)}/{2}\;}}\left( m \right)}{m!},$$
$${{q}_{0}}\left( x \right)=1,  \qquad{{q}_{2n}}\left( x \right)=\sum\limits_{m=1}^{n}{x{{\left( x+n+m-1 \right)}_{2m-1}}}\frac{{{s}_{n-m}}\left( 2m \right)}{\left( 2m \right)!},$$  
$${{q}_{2n+1}}\left( x \right)=\sum\limits_{m=0}^{n}{x{{\left( x+n+m \right)}_{2m}}}\frac{{{s}_{n-m}}\left( 2m+1 \right)}{\left( 2m+1 \right)!}.$$
{\bfseries Example 9.2.}
$$\left[ {{x}^{2n}} \right]{{R}^{\beta }}\left( x \right)=\sum\limits_{m=1}^{n}{\frac{\beta {{\left( \beta +n+m-1 \right)}_{2m-1}}}{\left( 2m \right)!}}\left( \begin{matrix}
   n+m-1  \\
   n-m  \\
\end{matrix} \right),$$
$$\left[ {{x}^{2n+1}} \right]{{R}^{\beta }}\left( x \right)=\sum\limits_{m=0}^{n}{\frac{\beta {{\left( \beta +n+m \right)}_{2m}}}{\left( 2m+1 \right)!}}\left( \begin{matrix}
   n+m  \\
   n-m  \\
\end{matrix} \right).$$
Denote
$${{u}_{0}}\left( \beta ,x \right)=1,    \qquad{{u}_{n}}\left( \beta ,x \right)=\sum\limits_{m=1}^{n}{\beta {{\left( \beta +\frac{n+m}{2}-1 \right)}_{m-1}}}\frac{{{s}_{{\left( n-m \right)}/{2}\;}}\left( m \right)}{m!}{{x}^{m}}.$$
Then
$${{\left( {{g}^{\left[ \varphi  \right]}}\left( x \right) \right)}^{\beta }}=\sum\limits_{n=0}^{\infty }{{{u}_{n}}}\left( \beta ,\varphi  \right){{x}^{n}}.$$

E-mail: {evgeniy\symbol{"5F}burlachenko@list.ru}
\end{document}